\documentclass[12pt,reqno]{amsart}

\usepackage[english]{babel}		
\usepackage[latin1]{inputenc}
\usepackage{amscd}
\usepackage{amsfonts}
\usepackage{amsgen}
\usepackage{amsmath}
\usepackage{amssymb}
\usepackage{amstext}
\usepackage{amsthm}
\usepackage{graphicx}
\usepackage{latexsym}
\usepackage{epsfig}
\usepackage{wrapfig}
\usepackage[all,knot,arc,import,poly]{xy}
\usepackage[usenames]{color}
\usepackage[mathscr]{eucal}
\usepackage{indentfirst}
\usepackage[T1]{fontenc}
\usepackage{setspace}
\usepackage{lscape}
\usepackage{tikz}
\usepackage{pgfplots}
\pgfplotsset{compat=1.13}
\usepackage{braids}
\usepackage[pdftex,                %
implicit=true,
bookmarks=true,
breaklinks=true,
bookmarksnumbered = true,
colorlinks= true,
linkcolor= blue,
urlcolor= green,
anchorcolor = yellow,
citecolor=blue,
]{hyperref}%
\everymath{\displaystyle}

\newtheorem{Def}{Definition}[section]
\newtheorem{Teo}[Def]{Theorem}

\newtheorem{Lema}[Def]{Lemma}
\newtheorem{Ex}[Def]{Example}
\newtheorem{Obs}[Def]{Remark}

\newtheorem{Cor}[Def]{Corollary}

\newcommand{\Z}{\mathbb{Z}}
\newcommand{\Aut}{\operatorname{Aut}}

\begin{document}

\title[conjugacy  problem and virtually cyclic subgroups in {$B_\MakeLowercase{n}/[P_\MakeLowercase{n},P_\MakeLowercase{n}]$}]{The conjugacy problem and virtually cyclic subgroups in the Artin braid group quotient $B_n/[P_n,P_n]$}
\author[O.~Ocampo]{Oscar Ocampo}
\address{Universidade Federal da Bahia, UFBA, Departamento de Matem\'{a}tica, 40170-110, Salvador, Bahia, Brazil}
\email{oscaro@ufba.br}

\author[P.~Santos Júnior]{Paulo Cesar Cerqueira dos Santos Júnior}
\address{Universidade Federal da Bahia, UFBA, Departamento de Matem\'{a}tica, 40170-110, Salvador, Bahia, Brazil}
\email{pcesarmath@gmail.com}

\date{\today}


\begin{abstract}

Let $n\geq 3$.  In this paper we deal with the conjugacy problem in the Artin braid group quotient $B_n/[P_n,P_n]$. To solve it we use systems of equations over the integers arising from the action of $B_n/[P_n,P_n]$ over the abelianization of the pure Artin braid group $P_n/[P_n,P_n]$. Using this technique we also realize explicitly infinite virtually cyclic subgroups in $B_n/[P_n,P_n]$.

\end{abstract}


\maketitle

\section{Introduction}











The Artin braid group $B_n$ was introduced by Artin in 1925 \cite{artin1} and further studied in 1947 \cite{artin2}, where a presentation for this group was given
$$
B_n= \left\langle \sigma_1,\ldots,\sigma_{n-1} \mid
\begin{array}{l}
\sigma_i\sigma_{i+1}\sigma_i=\sigma_{i+1}\sigma_i\sigma_{i+1} \\ \textrm{ for }1\leq i\leq n-2;\\
\sigma_i\sigma_j=\sigma_j\sigma_i \textrm{ for } \left|i-j\right|\geq 2
\end{array}
\right\rangle.
$$
The natural homomorphism $\sigma\colon B_n\longrightarrow S_n$, that relates a braid to a permutation, is defined on the given generators by $\sigma(\sigma_i)=(i, i+1)$ for all $1\le i\le n-1$. 
Since the operation in the braid group is geometrically defined as concatenation then, just as for braids, we read permutations from left to right so that if $\alpha,\beta\in S_n$ then their product is defined by $\alpha\cdot\beta(i)=\beta(\alpha(i))$ for $i=1,2,\dots,n$. The pure braid group $P_n$ is defined to be the kernel of $\sigma$. So, we obtain the following short exact sequence
\begin{equation}\label{spura}
1 \to  P_n \to  B_n\stackrel{\sigma}\to S_n \to 1.
\end{equation}
A generating set of $P_n$ is given by $\{A_{i,j}\}_{1\le i<j\le n}$, see \cite{artin1}, where
$$
A_{i,j}=\sigma_{j-1}\cdots\sigma_{i+1}\sigma_{i}^2\sigma_{i+1}^{-1}\cdots\sigma_{j-1}^{-1}.
$$

Let $[P_n,P_n]$ denote the commutator subgroup of $P_n$, recall that it is a characteristic subgroup of $P_n$. From equation~(\ref{spura}) we obtain the following short exact sequence
\begin{equation}\label{sqpura}
1 \to  P_n/[P_n,P_n] \to  B_n/[P_n,P_n]\stackrel{\overline{\sigma}} \to S_n \to 1,
\end{equation}  
where $\overline{\sigma}\colon B_n/[P_n,P_n]\longrightarrow S_n$ is the homomorphism induced by $\sigma\colon B_n\longrightarrow S_n$.
For $1\le i<j\le n$, we also set $A_{j,i}=A_{i,j}$, and if $A_{i,j}$ appears with exponent $m_{i,j}\in\Z$, then we let $m_{j,i}=m_{i,j}$. 
It was proved in \cite{art} that the Artin braid group quotient $B_n/[P_n,P_n]$ is a crystallographic group that has no 2-torsion elements and that its torsion is equal to the odd torsion of the symmetric group $S_n$. Furthermore, the authors studied the conjugacy classes of finite-order elements and showed that they are in correspondence with the conjugacy classes of elements of odd order in $S_n$, see \cite[Theorem~5]{art}. Hence, a remaining open question is about the conjugacy classes of any element of $B_n/[P_n,P_n]$, not necessarily of finite order.

\emph{The conjugacy problem} for a group $G$ is the following decision problem: to determine, given two words $x$ and $y$ in $G$, whether or not they represent conjugate elements of $G$. The conjugacy problem is one of the three famous decision problems in groups first formulated by Dehn in the early 20th century (see \cite{dehn1}).

The first solution for the conjugacy problem in braid groups was given by Garside in 1969 \cite{gar} and several successive improvements have been developed since, for instance Elrifai and Morton in 1994 \cite{em}, Birman, Ko and Lee in 1998 \cite{bkl1} and 2001 \cite{bkl2}, and Franco and Gonz\'alez-Meneses in 2003 \cite{fgm}. 
The common strategy to most solutions to the conjugacy problem in the braid group $B_n$, given two braids $x$ and $y$, is to compute finite conjugacy invariants and subsequently test whether they are disjoint sets. 
In \cite{m} the author exhibit seven different classes of braids in $B_3$ and proved that any element of $B_3$ is conjugate to a unique element of these seven classes. To classify conjugacy classes of elements in braid groups with several strings is not trivial, during the last years a lot of work has been done to create efficient algorithms to attack the problem of determine if two given braids are conjugate, see \cite{bkl1}, \cite{bkl2}, \cite{cw}, \cite{em} and \cite{fgm}. 


In the first part of this paper, Section~\ref{sec:2}, we study the conjugacy problem in the Artin braid group quotient $B_n/[P_n,P_n]$. The initial motivation comes from the fact that the conjugacy classes of finite-order elements of $B_n/[P_n,P_n]$ were studied in \cite[Theorem~5]{art}. In the crystallographic group $B_n/[P_n,P_n]$ we may solve the conjugacy problem using linear equations over the integers that arise from the action by conjugacy of $B_n/[P_n,P_n]$ over $P_n/[P_n,P_n]$. 


Let $n\geq 3$. To simplify the notation, an element $\overline{\beta}$ in $B_n/[P_n,P_n]$ will be denoted just by $\beta$. 
In what follows $T_{\theta}$ means a transversal of the action of $\langle \theta \rangle$ on the set $\mathcal{B}=\{A_{i,j}\mid 1\le i<j\le n\}$ and $\mathcal{O}_{\theta}(A_{r,s})$ means the orbit of the element $A_{r,s}$ under this action.
The main results of this paper about the conjugacy problem in the Artin braid group quotient $B_n/[P_n,P_n]$ are the following two theorems.

\begin{Teo}\label{teo1}
Let $n\ge 3$, let $\beta\in B_n/[P_n,P_n]$ and let $\theta=\overline{\sigma}(\beta^{-1})$. The elements $\beta\prod_{1\le i<j\le n}{A_{i,j}^{a_{i,j}}}$ and $\beta\prod_{A_{r,s}\in T_{\theta}}{A_{r,s}^{S_{r,s}}}$ are conjugate in $B_n/[P_n,P_n]$, where $a_{i,j}\in\Z$ for all $1\le i<j\le n$ and $T_{\theta}$ is a transversal of the action of $\langle\theta\rangle$ on the set $\mathcal{B}$ if and only if $S_{r,s}=\sum_{A_{i,j}\in\mathcal{O}_{\theta}(A_{r,s})}{a_{i,j}}$, for all $A_{r,s}\in T_{\theta}$.
\end{Teo}

Given any element $g$ of a group $G$ we denote by $C_G(g)$ the centralizer of $g$ in $G$.

\begin{Teo}\label{teo2}
Let $n\ge 3$, let $\beta\in B_n/[P_n,P_n]$ and let $\theta=\overline{\sigma}(\beta^{-1})$. 
Suppose that $\tilde{C}\in\overline{\sigma}^{-1}(C_{S_n}(\theta))$ and let $\overline{\sigma}(\tilde{C})=C$. 
Let $c_{i,j}\in\Z$, for $1\le i<j\le n$, such that $\tilde{C}\beta\tilde{C}^{-1}=\beta\prod_{1\le i<j\le n}{A_{i,j}^{c_{i,j}}}$. 
The elements $\beta\prod_{A_{t,q}\in T_{\theta}}{A_{t,q}^{z_{t,q}}}$ and $\beta\prod_{A_{t,q}\in T_{\theta}}{A_{t,q}^{b_{t,q}}}$, where $z_{t,q},b_{t,q}\in\Z$ for all $A_{t,q}\in T_{\theta}$, are conjugate in $B_n/[P_n,P_n]$ if and only if
\begin{enumerate}
	\item[$(i)$] $\sum_{A_{i,j}\in\mathcal{O}_{\theta}(A_{t,q})}{c_{i,j}}+z_{C(t),C(q)}=b_{t,q}$, if $A_{C^{-1}(t),C^{-1}(q)}\notin\mathcal{O}_{\theta}(A_{t,q})$ for $A_{t,q}\in T_{\theta}$, where $z_{C(t),C(q)}=\sum_{A_{i,j}\in\mathcal{O}_{\theta}(A_{C(t),C(q)})}{a_{i,j}}$ and $a_{i,j}\in\Z$ is the exponent of $A_{i,j}$;
	\item[$(ii)$] $\sum_{A_{i,j}\in\mathcal{O}_{\theta}(A_{t,q})}{c_{i,j}}+z_{t,q}=b_{t,q}$, if $A_{C^{-1}(t),C^{-1}(q)}\in\mathcal{O}_{\theta}(A_{t,q})$ for $A_{t,q}\in T_{\theta}$.
\end{enumerate}
\end{Teo}


Now, in a second part of this work we deal with a problem of realization of groups in $B_n/[P_n,P_n]$. In \cite{art} was proved that every finite abelian subgroup of odd order of $S_n$ may be realized in $B_n/[P_n,P_n]$, and also the authors proved that for $n\geq 7$ the non-abelian group of order 21 (that belongs to a family called Frobenius group) is a subgroup of $B_n/[P_n,P_n]$. The realization problem of subgroups is interesting by itself, the paper \cite{art2} is devoted entirely to realize finite subgroups in quotients of the Artin braid group, in particular in the group $B_n/[P_n,P_n]$.


In Section~\ref{sec:3}, we are interested in to realize in $B_n/[P_n,P_n]$ a family of (infinite) subgroups called virtually cyclic groups, that are groups with a cyclic subgroup $H$ of finite index. Virtually cyclic groups are a very interesting family of groups; for instance, every abelian subgroup of a Gromov hyperbolic group is virtually cyclic (see \cite{a}), and an infinite group is virtually cyclic if and only if it is finitely generated and has exactly two ends (see \cite{e}).
The classification of virtually cyclic subgroups in a given group $G$ is an interesting problem in its own right. Besides helping us to understand better the structure of $G$, in the case $G$ satisfies the Fibered Isomorphism Conjecture of Farrell and Jones (see \cite{fajo}), the algebraic $K$-theory of their group rings (over $\Z$) may be computed by means of the algebraic $K$-theory groups of their virtually cyclic subgroups, via the so-called ``assembly maps''. For more details about this topic we refer the reader to \cite{virtual}, \cite{gg} and \cite{gjm} and their references. 
We note that the classification of virtually cyclic subgroups of surface braid groups was studied for the case of the sphere in  \cite{virtual} and for the projective plane in \cite{gg}. The lower algebraic $K$-theory of the sphere braid groups was computed in \cite{gjm}.   

Since the group $B_n/[P_n,P_n]$ has no 2-torsion (\cite[Theorem~2]{art}) then the only (infinite) virtually cyclic subgroups of this quotient are of the form semi-direct product, because the amalgamated products have even order elements, see Theorem~\ref{teo17}. In Theorem~\ref{teoprimos} we realize virtually cyclic subgroups in $B_n/[P_n,P_n]$ of the form $\Z_p\rtimes\Z$, where $p\leq n$ is an odd prime number.

\begin{Teo}\label{teoprimos}
Let $n\geq 3$. Let $p\leq n$ be an odd prime number. The subgroup $\Z_p\rtimes_{\alpha_k}\Z$ can be realized in $B_n/[P_n,P_n]$, where $k=1,\dots,p-1$. 
\end{Teo}

\subsection*{Acknowledgments}
The authors would like to thank Daciberg Lima Gon\c calves for helpful discussions and for suggest us the realization problem of virtually cyclic groups. The second author was partially supported by CAPES.

\section{The conjugacy problem and conjugacy classes of elements in $B_n/[P_n,P_n]$}\label{sec:2}

In this section we shall prove Theorems~\ref{teo1} and~\ref{teo2},  in order to determine conditions for two elements of $B_n/[P_n,P_n]$ being conjugate.


\begin{Obs}\label{rmk1}
Let $\beta\in B_n/[P_n,P_n]$. A necessary condition for a given coset of a braid $y\in B_n/[P_n,P_n]$ to belong to the conjugacy class of $\beta$ is that the permutations $\overline{\sigma}(y)$ and $\overline{\sigma}(\beta)$ need to have the same cycle type, where $\overline{\sigma}$ is the surjective homomorphism given in (\ref{sqpura}).
Now, consider $y\in B_n/[P_n,P_n]$ such that $\overline{\sigma}(y)$ and $\overline{\sigma}(\beta)$ have the same cycle type. So, there is $\hat{\gamma}\in S_n$ such that $\hat{\gamma}\overline{\sigma}(y)\hat{\gamma}^{-1}=\overline{\sigma}(\beta)$. 
Since $\overline{\sigma}$ is surjective, then there is $\gamma \in B_n/[P_n,P_n]$ such that $\tilde{y}=\gamma y \gamma^{-1}$ belongs to $\overline{\sigma}^{-1}(\overline{\sigma}(\beta))$.
We note that $\tilde{y}$ is not necessarily conjugate to $\beta$.
For instance, let $\beta=\sigma_2\sigma_1^{-1}$ and $y=\sigma_1\sigma_2$ in $B_3/[P_3,P_3]$. Notice that $\overline{\sigma}(\beta)=(1,2,3)$ and $\overline{\sigma}(y)=(3,2,1)$. Consider $\hat{\gamma}=(1,3)$ in $S_3$ and $\gamma=\sigma_1\sigma_2\sigma_1^{-1}$ in $B_3/[P_3,P_3]$ so, $\hat{\gamma}\overline{\sigma}(y){\hat{\gamma}}^{-1}=\overline{\sigma}(\beta)$ and $\overline{\sigma}(\gamma)=\hat{\gamma}$. Thus, $\tilde{y}=\gamma y \gamma^{-1}$ belongs to $\overline{\sigma}^{-1}(\overline{\sigma}(\beta))$. However $\tilde{y}$ is not conjugate to $\beta$, since $\beta$ has finite order (see \cite[Lemma~28]{art}) while $\tilde{y}$ is not a finite order element in $B_3/[P_3,P_3]$ (see \cite[Proof of Proposition~21]{art}).
\end{Obs}

From now onwards in this section we deal with the problem to decide whether two cosets of braids of the Artin braid group quotient $B_n/[P_n,P_n]$, with the same permutation, are conjugate.


The action by conjugation of $B_n/[P_n,P_n]$ on $P_n/[P_n,P_n]$ is given by 
\begin{equation}\label{eq:1.1}
\alpha A_{i,j}\alpha^{-1}=A_{\pi(i),\pi(j)},
\end{equation}
where $\alpha\in B_n/[P_n,P_n]$ and $\pi=\overline{\sigma}(\alpha^{-1})$, see~\cite[Proposition~12]{art}. Furthermore by \cite[Remark~16]{art} given $\alpha\in B_n/[P_n,P_n]$ and $\prod_{1\le i< j\le n}{A_{i,j}^{m_{i,j}}}\in P_n/[P_n,P_n]$ we have 
\begin{equation}\label{eq:1.2}
\alpha\left(\prod_{1\le i< j\le n}{A_{i,j}^{m_{i,j}}}\right)\alpha^{-1}=\prod_{1\le i< j\le n}{A_{\pi(i),\pi(j)}^{m_{i,j}}}=\prod_{1\le i< j\le n}{A_{i,j}^{m_{\pi^{-1}(i),\pi^{-1}(j)}}},
\end{equation} 
where $m_{i,j}\in\Z$ and $\pi=\overline{\sigma}(\alpha^{-1})$.

Let $\beta\in B_n/[P_n,P_n]$ such that $\overline{\sigma}(\beta^{-1})=\theta\in S_n$. The action by conjugation of $B_n/[P_n,P_n]$ on $P_n/[P_n,P_n]$ induce an action $\langle\theta\rangle\times\mathcal{B}\longrightarrow \mathcal{B}$ where $\mathcal{B}=\{A_{i,j}\mid 1\le i<j\le n\}$ and is given by 
\begin{equation}\label{eq:1.3}
\theta^k(A_{i,j})=A_{\theta^k(i),\theta^k(j)},
\end{equation}
for $1\le k\le |\theta|$, where $|\theta|$ is the order of $\theta$. We denote by $\mathcal{O}_{\theta}(A_{i,j})$ the orbit of the element $A_{i,j}$ and by $T_\theta$ a transversal for this action. In what follows given a set $X$ we denote its cardinality by $|X|$. Now, we prove Theorem~\ref{teo1}.




\begin{proof}[Proof of Theorem~\ref{teo1}]
Let $n\ge 3$, let $\beta\in B_n/[P_n,P_n]$ and let $\theta=\overline{\sigma}(\beta^{-1})$. 
If there is $X=\prod_{1\le i<j\le n}{A_{i,j}^{x_{i,j}}}$ in $P_n/[P_n,P_n]$, where $x_{i,j}\in \Z$ are indeterminates, such that $X\beta\prod_{1\le i<j\le n}{A_{i,j}^{a_{i,j}}}X^{-1}=\beta\prod_{A_{r,s}\in T_{\theta}}{A_{r,s}^{S_{r,s}}}$ then the elements $\beta\prod_{1\le i<j\le n}{A_{i,j}^{a_{i,j}}}$ and $\beta\prod_{A_{r,s}\in T_{\theta}}{A_{r,s}^{S_{r,s}}}$ are conjugate in $B_n/[P_n,P_n]$. 
Notice that, using the action described in equation~\eqref{eq:1.1}, we have
$$
\begin{array}{l}
X\beta\prod_{1\le i<j\le n}{A_{i,j}^{a_{i,j}}}X^{-1}=\beta\prod_{A_{r,s}\in T_{\theta}}{A_{r,s}^{S_{r,s}}}\iff X\beta\left(\prod_{1\le i<j\le n}{A_{i,j}^{a_{i,j}}}\right)X^{-1}\left(\prod_{A_{r,s}\in T_{\theta}}{A_{r,s}^{-S_{r,s}}}\right)\beta^{-1}=1\\
\iff X\left(\prod_{1\le i<j\le n}{A_{\theta(i),\theta(j)}^{a_{i,j}}}\right)\beta X^{-1}\beta^{-1}\prod_{A_{r,s}\in T_{\theta}}{A_{\theta(r),\theta(s)}^{-S_{r,s}}}=1\\
\iff \prod_{1\le i<j\le n}{A_{i,j}^{x_{i,j}}}\prod_{1\le i<j\le n}{A_{\theta(i),\theta(j)}^{a_{i,j}}} 
\prod_{1\le i<j\le n}{A_{i,j}^{-x_{\theta^{-1}(i),\theta^{-1}(j)}}}\prod_{A_{r,s}\in T_{\theta}}{A_{\theta(r),\theta(s)}^{-S_{r,s}}}=1\\
\iff \prod_{1\le i<j\le n}{A_{i,j}^{x_{i,j}}}\prod_{1\le i<j\le n}{A_{i,j}^{a_{\theta^{-1}(i),\theta^{-1}(j)}}} 
\prod_{1\le i<j\le n}{A_{i,j}^{-x_{\theta^{-1}(i),\theta^{-1}(j)}}}\prod_{A_{r,s}\in T_{\theta}}{A_{\theta(r),\theta(s)}^{-S_{r,s}}}=1.
\end{array}
$$
Since $P_n/[P_n,P_n]$ is a torsion free abelian group with basis $\{A_{i,j}\}_{1\leq i<j\leq n}$, the integer solutions of the last equation are the ones of the following system of equations
\begin{equation}\label{sisgeral}
\left\{\begin{array}{lc}
x_{i,j}+a_{\theta^{-1}(i),\theta^{-1}(j)}-x_{\theta^{-1}(i),\theta^{-1}(j)}=0\\
x_{\theta(r),\theta(s)}+a_{r,s}-x_{r,s}-S_{r,s}=0.
\end{array}\right.
\end{equation}
for all $1\le i<j\le n$ and $i\neq\theta(r)$, $j\neq\theta(s)$ and for all $\{r,s\}$ such that $A_{r,s}\in T_{\theta}$. 
Note that system~(\ref{sisgeral}) have solution if and only if the following subsystems~(\ref{sisorb}) have solution
\begin{equation}\label{sisorb}
\left\{\begin{array}{lc}
x_{\theta^t(r),\theta^t(s)}+a_{\theta^{-(1+m-t)}(r),\theta^{-(1+m-t)}(s)}-x_{\theta^{-(1+m-t)}(r),\theta^{-(1+m-t)}(s)}=0\\
x_{\theta(r),\theta(s)}+a_{r,s}-x_{r,s}-S_{r,s}=0.
\end{array}\right.
\end{equation}
for all $2\le t\le m$ and for all $\{r,s\}$ such that $A_{r,s}\in T_{\theta}$ where $m=|\mathcal{O}_{\theta}(A_{r,s})|$ 

Let $A_{r,s}\in T_{\theta}$ and consider the elements $e_1=A_{r,s},e_2=A_{\theta(r),\theta(s)},\dots,e_m=A_{\theta^{m-1}(r),\theta^{m-1}(s)}$ where $m=|\mathcal{O}_{\theta}(A_{r,s})|$. 
Thus the permutation matrix of the orbit of $A_{r,s}$ by the action of $\langle\theta\rangle$ with respect to the ordered set $\{e_1,e_2,\dots,e_m\}$ is given by the matrix $M_{r,s}$ 
$$M_{r,s}=
\begin{pmatrix}
0&0&\ldots&0&1\\
1&0&\ldots&0&0\\
0&1&\ldots&0&0\\
\vdots&\vdots&\ddots&\vdots&\vdots\\ 
0&0&\ldots&1&0
\end{pmatrix}.
$$
So, for all $\{r,s\}$ such that $A_{r,s}\in T_{\theta}$ and considering $m=|\mathcal{O}_{\theta}(A_{r,s})|$, we can write the subsystems of equations~(\ref{sisorb}) as follows 
$$
(\star)\quad [I_{r,s}-M_{r,s}]
\begin{pmatrix}
x_{r,s}\\
x_{\theta(r),\theta(s)}\\
\vdots\\
x_{\theta^{m-1}(r),\theta^{m-1}(s)}
\end{pmatrix}=
\begin{pmatrix}
-a_{\theta^{-1}(r),\theta^{-1}(s)}\\
-a_{r,s}+S_{r,s}\\
\vdots\\
-a_{\theta^{m-2}(r),\theta^{m-2}(s)}
\end{pmatrix}
$$
where $I_{r,s}=(Id)_{m\times m}$. Analyzing the coefficient matrix and the extended matrix associated with this system, we concluded that $(\star)$ have solution if and only if $S_{r,s}=\sum_{A_{i,j}\in\mathcal{O}_{\theta}(A_{r,s})}{a_{i,j}}$. Thus the elements $\beta\prod_{1\le i<j\le n}{A_{i,j}^{a_{i,j}}}$ and $\beta\prod_{A_{r,s}\in T_{\theta}}{A_{r,s}^{S_{r,s}}}$ are conjugate in $B_n/[P_n,P_n]$ if and only if $S_{r,s}=\sum_{A_{i,j}\in\mathcal{O}_{\theta}(A_{r,s})}{a_{i,j}}$. 
\end{proof}

\begin{Obs}

With the notations of Theorem~\ref{teo1} let $T'_{\theta}$ be another transversal of the action of $\langle\theta\rangle$ on the set $\mathcal{B}$. Supppose that $A_{r',s'}\in T'_{\theta}$ is the representative of the orbit of $A_{r,s}\in T_{\theta}$. From Theorem~\ref{teo1}, if $S_{r,s}=\sum_{A_{i,j}\in\mathcal{O}_{\theta}(A_{r,s})}{a_{i,j}}=\sum_{A_{i,j}\in\mathcal{O}_{\theta}(A_{r',s'})}{a_{i,j}}$ then $\beta\prod_{A_{r,s}\in T_{\theta}}{A_{r,s}^{S_{r,s}}}$ and $\beta\prod_{A_{r,s}\in T_{\theta}}{A_{r',s'}^{S_{r,s}}}$ are conjugate.

\end{Obs}

Let $\beta\in B_n/[P_n,P_n]$ satisfying the conditions of Theorem~\ref{teo1}. Let $\overline{\sigma}(\beta^{-1})=\theta$, we choose a transversal $T_{\theta}$ of the action of $\langle \theta \rangle$ on the set $\mathcal{B}=\{A_{i,j}\mid 1\le i<j\le n\}$. 
We would like to know which kind of elements belongs to the conjugacy class of $\beta\prod_{A_{r,s}\in T_{\theta}}{A_{r,s}^{S_{r,s}}}$, having the same description using the chosen transversal $T_{\theta}$, where $S_{r,s}=\sum_{A_{i,j}\in\mathcal{O}_{\theta}(A_{r,s})}{a_{i,j}}$ for all $A_{r,s}\in T_{\theta}$. 
In the next two examples we shall illustrate the procedure to finding conjugate elements to the coset of the braid $\beta\prod_{A_{r,s}\in T_{\theta}}{A_{r,s}^{z_{r,s}}}$, given $\beta\in B_n/[P_n,P_n]$ and $T_{\theta}$ fixed and $z_{r,s}\in \Z$. 
This method is as follows. 
Suppose that $\beta\prod_{A_{r,s}\in T_{\theta}}{A_{r,s}^{z_{r,s}}}$ and $\beta\prod_{A_{r,s}\in T_{\theta}}{A_{r,s}^{b_{r,s}}}$ are two arbitrary elements in $\overline{\sigma}^{-1}(\overline{\sigma}(\beta))$. 
Let $X=\prod_{1\le i<j\le n}{A_{i,j}^{x_{i,j}}}$, with $x_{i,j}\in \Z$ fol all $1\leq i<j\leq n$. 
Consider $\tilde{C}\in\overline{\sigma}^{-1}(C_{S_n}(\theta))$.
We highlight that the requirement to the element $\tilde{C}$ belongs to $\overline{\sigma}^{-1}(C_{S_n}(\theta))$ comes from the fact that $\overline{\sigma}(X\tilde{C}\beta\prod_{A_{r,s}\in T_{\theta}}{A_{r,s}^{z_{r,s}}}\tilde{C}^{-1}X^{-1})$ need to be equal to $\overline{\sigma}(\beta\prod_{A_{r,s}\in T_{\theta}}{A_{r,s}^{b_{r,s}}})$.
Hence we must find conditions for $X$ such that the equation $X\tilde{C}\beta\prod_{A_{r,s}\in T_{\theta}}{A_{r,s}^{z_{r,s}}}\tilde{C}^{-1}X^{-1}=\beta\prod_{A_{r,s}\in T_{\theta}}{A_{r,s}^{b_{r,s}}}$ is true. 

\begin{Ex}\label{exm}
Let $\beta=\sigma_1\sigma_2\in B_6/[P_6,P_6]$ and consider $\overline{\sigma}((\sigma_1\sigma_2)^{-1})=\theta=(1,2,3)$. Note that the element $\sigma_4\sigma_5$ is an element of the centralizer of $\sigma_1\sigma_2$ in $B_6/[P_6,P_6]$. Fix the transversal 
 $T_{\theta}=\{A_{1,2},A_{1,4},A_{1,5},A_{1,6},A_{4,5},A_{4,6},A_{5,6}\}$ of the action of $\langle\theta\rangle$ on the set $\mathcal{B}=\{A_{i,j}\mid 1\le i<j\le 6\}$. 
Let $X=\prod_{1\le i<j\le 6}{A_{i,j}^{x_{i,j}}}$ and let $\overline{\sigma}(\sigma_4\sigma_5)=C$. So, applying the action described in equations~\eqref{eq:1.1} and~\eqref{eq:1.2}, we have 
\begin{align*}
& X(\sigma_4\sigma_5)\sigma_1\sigma_2\left(\prod_{A_{t,q}\in T_{\theta}}{A_{t,q}^{z_{t,q}}}\right)(\sigma_4\sigma_5)^{-1}X^{-1}  
=\sigma_1\sigma_2\prod_{A_{t,q}\in T_{\theta}}{A_{t,q}^{b_{t,q}}} \\
 \iff & X(\sigma_4\sigma_5)\sigma_1\sigma_2\left(\prod_{A_{t,q}\in T_{\theta}}{A_{t,q}^{z_{t,q}}}\right)(\sigma_4\sigma_5)^{-1} 
X^{-1}\left(\prod_{A_{t,q}\in T_{\theta}}{A_{t,q}^{-b_{t,q}}}\right)(\sigma_1\sigma_2)^{-1}=1\\
\iff& X(\sigma_4\sigma_5)\sigma_1\sigma_2(\sigma_4\sigma_5)^{-1}\prod_{A_{t,q}\in T_{\theta}}{A_{C^{-1}(t),C^{-1}(q)}^{z_{t,q}}}
  X^{-1}\left(\prod_{A_{t,q}\in T_{\theta}}{A_{t,q}^{-b_{t,q}}}\right)(\sigma_1\sigma_2)^{-1}=1\\
\iff& X\sigma_1\sigma_2\left(\prod_{A_{t,q}\in T_{\theta}}{A_{C^{-1}(t),C^{-1}(q)}^{z_{t,q}}}\right)X^{-1}
  \left(\prod_{A_{t,q}\in T_{\theta}}{A_{t,q}^{-b_{t,q}}}\right)(\sigma_1\sigma_2)^{-1}=1\\
\iff& X\left(\prod_{A_{t,q}\in T_{\theta}}{A_{\theta(C^{-1}(t)),\theta(C^{-1}(q))}^{z_{t,q}}}\right)\sigma_1\sigma_2X^{-1}
 (\sigma_1\sigma_2)^{-1}\prod_{A_{t,q}\in T_{\theta}}{A_{\theta(t),\theta(q)}^{-b_{t,q}}}=1\\
\iff& X\left(\prod_{A_{t,q}\in T_{\theta}}{A_{\theta(C^{-1}(t)),\theta(C^{-1}(q))}^{z_{t,q}}}\right)\prod_{1\le i<j\le 6}{A_{i,j}^{-x_{\theta^{-1}(i),\theta^{-1}(j)}}}
 \prod_{A_{t,q}\in T_{\theta}}{A_{\theta(t),\theta(q)}^{-b_{t,q}}}=1.
\end{align*}
Since $P_n/[P_n,P_n]$ is a torsion free abelian group with basis $\{A_{i,j}\}_{1\leq i<j\leq n}$, then the last equation is true if the following systems are true, where each system is relative to the exponent of each $A_{t,q}\in T_{\theta}$
\begin{eqnarray}\label{tip1}
\left\{\begin{array}{lc}
x_{1,2}-x_{1,3}=0\\
x_{2,3}+z_{1,2}-b_{1,2}-x_{1,2}=0\\
x_{1,3}-x_{2,3}=0
\end{array}\right.
\end{eqnarray}
\begin{eqnarray}\label{tip2}
\left\{\begin{array}{ll}
x_{1,4}-x_{3,4}=0\\
x_{2,4}-b_{1,4}-x_{1,4}+z_{1,6}=0\\
x_{3,4}-x_{2,4}=0
\end{array}\right.,\,
\left\{\begin{array}{ll}
x_{1,5}-x_{3,5}=0\\
x_{2,5}-b_{1,5}-x_{1,5}+z_{1,4}=0\\
x_{3,5}-x_{2,5}=0
\end{array}\right.,\,
\left\{\begin{array}{ll}
x_{1,6}-x_{3,6}=0\\
x_{2,6}-b_{1,6}-x_{1,6}+z_{1,5}=0\\
x_{3,6}-x_{2,6}=0
\end{array}\right.
\end{eqnarray}
\begin{eqnarray}\label{tip3}
\left\{\begin{array}{lc}
x_{4,5}-x_{4,5}+z_{4,6}-b_{4,5}=0
\end{array}\right.,\,
\left\{\begin{array}{lc}
x_{4,6}-x_{4,6}+z_{5,6}-b_{4,6}=0
\end{array}\right.
\end{eqnarray}
\begin{eqnarray}
\left\{\begin{array}{lc}\nonumber
x_{5,6}-x_{5,6}+z_{4,5}-b_{5,6}=0.
\end{array}\right.
\end{eqnarray}

The system of equations~(\ref{tip1}) have solution if and only if $z_{1,2}=b_{1,2}$. For the systems of equations given in~(\ref{tip2}) to have solution, it is necessary that $b_{t,q}=z_{C(t),C(q)}$ where $(t,q)\in\{(1,4),(1,5),(1,6)\}$. The system~(\ref{tip3}) is trivially solved, just consider $z_{4,6}=b_{4,5}$, $z_{5,6}=b_{4,6}$ and $z_{4,5}=b_{5,6}$. 
Therefore, we found a condition for elements of the form $\beta\prod_{A_{t,q}\in T_{\theta}}{A_{t,q}^{b_{t,q}}}$ to be in the conjugacy class of $\beta\prod_{A_{t,q}\in T_{\theta}}{A_{t,q}^{z_{t,q}}}$. 
For instance, the element $\sigma_1\sigma_2A_{1,2}^{z_{1,2}}A_{1,4}^{z_{1,6}}A_{1,5}^{z_{1,4}}A_{1,6}^{z_{1,5}}A_{4,5}^{z_{4,6}}A_{4,6}^{z_{5,6}}A_{5,6}^{z_{4,5}}$ belongs to the conjugacy class of $\sigma_1\sigma_2\prod_{A_{t,q}\in T_{\theta}}{A_{t,q}^{z_{t,q}}}$. 

\end{Ex}

In Example~\ref{exm} we considered an element of the centralizer of $\sigma_1\sigma_2$ in $B_6/[P_6,P_6]$, but in general we must consider an element in $\overline{\sigma}^{-1}(C_{S_n}(\theta))$ instead $C_{B_n/[P_n,P_n]}(\theta)$. We will introduce this situation, more generally, in the next example.

\begin{Ex}

Let $\beta=\sigma_1\sigma_2\in B_6/[P_6,P_6]$ and consider $\theta=\overline{\sigma}((\sigma_1\sigma_2)^{-1})=(1,2,3)$. Note that $\sigma_4A_{1,2}\in\overline{\sigma}^{-1}(C_{S_6}(\theta))$. Fix the transversal $T_{\theta}=\{A_{1,2},A_{1,4},A_{1,5},A_{1,6},A_{4,5},A_{4,6},A_{5,6}\}$ of the action of $\langle\theta\rangle$ in the set $\mathcal{B}$.
Doing a similar computation as the one given in Example~\ref{exm}, we found a condition for elements of the form $\beta\prod_{A_{t,q}\in T_{\theta}}{A_{t,q}^{b_{t,q}}}$ to be in the conjugacy class of $\beta\prod_{A_{t,q}\in T_{\theta}}{A_{t,q}^{z_{t,q}}}$. 
In fact, necessary conditions are $z_{1,2}=b_{1,2}$ and $b_{t,q}=z_{C(t),C(q)}$, where $\overline{\sigma}(\sigma_4A_{1,2})=C$ and $(t,q)\in\{(1,4),(1,5),(1,6)\}$.  
For instance, the element $\sigma_1\sigma_2A_{1,2}^{z_{1,2}}A_{1,4}^{z_{1,5}}A_{1,5}^{z_{1,4}}A_{1,6}^{z_{1,6}}A_{4,5}^{z_{4,5}}A_{4,6}^{z_{5,6}}A_{5,6}^{z_{4,6}}$ is in the conjugacy class of $\sigma_1\sigma_2\prod_{A_{t,q}\in T_{\theta}}{A_{t,q}^{z_{t,q}}}$.

\end{Ex}


Although the next result is simple it will be useful in the proof of Theorem~\ref{teo2}.

\begin{Lema}\label{lemma2.6}
Let $\beta\in B_n/[P_n,P_n]$ and let $\theta=\overline{\sigma}(\beta^{-1})$. Given $\tilde{C}\in\overline{\sigma}^{-1}(C_{S_n}(\theta))$, there are $c_{i,j}\in\Z$ such that $\tilde{C}\beta\tilde{C}^{-1}=\beta\prod_{1\le i<j\le n}{A_{i,j}^{c_{i,j}}}$.
\end{Lema}

\begin{proof}
By hypothesis follows that $\overline{\sigma}(\tilde{C}\beta\tilde{C}^{-1})=\theta$, so $\overline{\sigma}(\beta^{-1}\tilde{C}\beta\tilde{C}^{-1})=1$. Thus, there are $c_{i,j}\in\Z$ such that $\beta^{-1}\tilde{C}\beta\tilde{C}^{-1}=\prod_{1\le i<j\le n}{A_{i,j}^{c_{i,j}}}$. 
\end{proof}

The information of Lemma~\ref{lemma2.6} was used in the statement of Theorem~\ref{teo2}, that we now prove. 

\begin{proof}[Proof of Theorem~\ref{teo2}]
Let $X=\prod_{1\le i<j\le n}{A_{i,j}^{x_{i,j}}}$. Suppose that $\tilde{C}\in \overline{\sigma}^{-1}(C_{S_n}(\theta))$ and let $\overline{\sigma}(\tilde{C})=C$. Hence, applying equation~\eqref{eq:1.2} and Lemma~\ref{lemma2.6}, we have
\begin{align*}
& X\tilde{C}\beta\left(\prod_{A_{t,q}\in T_{\theta}}{A_{t,q}^{z_{t,q}}}\right)\tilde{C}^{-1}X^{-1}=\beta\prod_{A_{t,q}\in T_{\theta}}{A_{t,q}^{b_{t,q}}}\\
\iff& X\tilde{C}\beta\left(\prod_{A_{t,q}\in T_{\theta}}{A_{t,q}^{z_{t,q}}}\right)\tilde{C}^{-1}X^{-1}
\left(\prod_{A_{t,q}\in T_{\theta}}{A_{t,q}^{-b_{t,q}}}\right)\beta^{-1}=1\\
\iff&X\tilde{C}\beta\tilde{C}^{-1}\left(\prod_{A_{t,q}\in T_{\theta}}{A_{t,q}^{z_{C(t),C(q)}}}\right)X^{-1}
\left(\prod_{A_{t,q}\in T_{\theta}}{A_{t,q}^{-b_{t,q}}}\right)\beta^{-1}=1\\
\end{align*}
\begin{align*}
\iff&X\beta\prod_{1\le i<j\le n}{A_{i,j}^{c_{i,j}}}\left(\prod_{A_{t,q}\in T_{\theta}}{A_{t,q}^{z_{C(t),C(q)}}}\right)X^{-1}
\left(\prod_{A_{t,q}\in T_{\theta}}{A_{t,q}^{-b_{t,q}}}\right)\beta^{-1}=1\\
\iff & X\prod_{1\le i<j\le n}{A_{\theta(i),\theta(j)}^{c_{i,j}}}\left(\prod_{A_{t,q}\in T_{\theta}}{A_{\theta(t),\theta(q)}^{z_{C(t),C(q)}}}\right)
\beta X^{-1}\beta^{-1}\prod_{A_{t,q}\in T_{\theta}}{A_{\theta(t),\theta(q)}^{-b_{t,q}}}=1\\
\iff&X\prod_{1\le i<j\le n}{A_{i,j}^{c_{\theta^{-1}(i),\theta^{-1}(j)}}}\left(\prod_{A_{t,q}\in T_{\theta}}{A_{\theta(t),\theta(q)}^{z_{C(t),C(q)}}}\right)
\beta X^{-1}\beta^{-1}\prod_{A_{t,q}\in T_{\theta}}{A_{\theta(t),\theta(q)}^{-b_{t,q}}}=1.
\end{align*}

Suppose that $A_{C^{-1}(t),C^{-1}(q)}\notin\mathcal{O}_{\theta}(A_{t,q})$ for $A_{t,q}\in T_{\theta}$. First we prove that $z_{C(t),C(q)}$, the exponent of $A_{C(t),C(q)}$, is equal to the exponent of an element of the transversal $T_{\theta}$. In case $A_{C(t),C(q)}\in T_{\theta}$ it is nothing to prove. 
Now suppose that $A_{C(t),C(q)}\notin T_{\theta}$, then there exists $A_{t_1,q_1}\in T_{\theta}$ such that $A_{C(t),C(q)}\in\mathcal{O}_{\theta}(A_{t_1,q_1})$, because $T_{\theta}$ is a transversal of the action of $\langle \theta \rangle$ on the set $\mathcal{B}=\{A_{i,j}\mid 1\le i<j\le n\}$. 
From the conditions given in the first item of this theorem we are assuming that $z_{C(t),C(q)}=\sum_{A_{i,j}\in\mathcal{O}_{\theta}(A_{C(t),C(q)})}{a_{i,j}}$, where $a_{i,j}\in\Z$ is the exponent of $A_{i,j}$. Since $A_{C(t),C(q)}\in\mathcal{O}_{\theta}(A_{t_1,q_1})$ then we obtain $z_{C(t),C(q)}=z_{t_1,q_1}$. 
Notice that we have the action of $\langle C\rangle$ on the set $\mathcal{B}$ and $\langle C\rangle\cdot T_{\theta}=T_{1_{\theta}}\cup T_{2_{\theta}}$, where  $T_{1_{\theta}}$ are the representatives of the orbits that are no invariant by $\langle C\rangle$ and $T_{2_{\theta}}$ are the representatives of the orbits preserved by $\langle C\rangle$. 
Let $T_{3_{\theta}}$ the set obtained by the action of $\langle C\rangle$ on $T_{1_{\theta}}$, $T_{3_{\theta}}=\langle C\rangle\cdot T_{1_{\theta}}$
and is not difficult to prove that $T'_{\theta}=T_{3_{\theta}}\cup T_{2_{\theta}}$ is also a transversal of the action of $\langle \theta\rangle$ on the set $\mathcal{B}$.
So the last equation is equivalent to the following system of equations 
\begin{align}
&\left\{\begin{array}{lc}\label{eqteo2}
x_{\theta^k(t),\theta^k(q)}-x_{\theta^{-(1+m-k)}(t),\theta^{-(1+m-k)}(q)}+c_{\theta^{-(1+m-k)}(t),\theta^{-(1+m-k)}(q)}=0\\
x_{\theta(t),\theta(q)}+c_{t,q}+z_{C(t),C(q)}-b_{t,q}-x_{t,q}=0
\end{array}\right.
\end{align}
for all $2\le k\le m$ and for all $A_{t,q}\in T_{\theta}$ where $m=|\mathcal{O}_{\theta}(A_{t,q})|$. 
Following the argument given in the proof of Theorem~\ref{teo1} we may conclude that equation~(\ref{eqteo2}) is equivalent to
$$
(I)\quad [I_{t,q}-M_{t,q}]
\begin{pmatrix}
x_{t,q}\\
x_{\theta(t),\theta(q)}\\
\vdots\\
x_{\theta^{m-1}(t),\theta^{m-1}(q)}
\end{pmatrix}=
\begin{pmatrix}
-c_{\theta^{-1}(t),\theta^{-1}(q)}\\
-c_{t,q}-z_{C(t),C(q)}+b_{t,q}\\
\vdots\\
-c_{\theta^{m-2}(t),\theta^{m-2}(q)}
\end{pmatrix}
$$
where $m=|\mathcal{O}_{\theta}(A_{t,q})|$, $I_{t,q}=(Id)_{m\times m}$ and $M_{t,q}$ is the permutation matrix of the orbit of $A_{t,q}$ by the action of $\langle\theta\rangle$ in the set $\mathcal{B}$ with respect to the ordered set $\{e_1=A_{t,q},e_2=A_{\theta(t),\theta(q)},\dots,e_m=A_{\theta^{m-1}(t),\theta^{m-1}(q)}\}$. Analyzing the coefficient matrix and the extended matrix associated to this system, we concluded that $(I)$ have solution if and only if $\sum_{A_{i,j}\in\mathcal{O}_{\theta}(A_{t,q})}{c_{i,j}}+z_{C(t),C(q)}=b_{t,q}$, proving item $(i)$. 
Similarly, if $A_{C^{-1}(t),C^{-1}(q)}\in\mathcal{O}_{\theta}(A_{t,q})$ for $A_{t,q}\in T_{\theta}$ the system $(I)$ have solution if and only if $\sum_{A_{i,j}\in\mathcal{O}_{\theta}(A_{t,q})}{c_{i,j}}+z_{t,q}=b_{t,q}$, verifying item $(ii)$.
\end{proof}





Next, we describe the steps to compare whether two given braids are conjugate in $B_n/[P_n,P_n]$. Let $\beta\in B_n/[P_n,P_n]$. Let $\alpha$ be an arbitrary braid different from $\beta$, we would like to decide whether $\alpha$ is conjugate to $\beta$ in $B_n/[P_n,P_n]$. Follows from the existence of the short exact sequence given in equation~(\ref{sqpura}) that a necessary condition to $\alpha$ and $\beta$ to be conjugate is that their permutations related have the same cycle type.
%
\begin{itemize}
	\item[Step 1.] Since $\alpha$ and $\beta$ have permutations with the same cycle type and $\overline{\sigma}$ is surjective 
	then $\alpha$ is conjugate to a braid $\alpha_1$ in $\overline{\sigma}^{-1}(\overline{\sigma}(\beta))$, see Remark~\ref{rmk1}.
\item[Step 2.]  
Suppose that $\alpha_1=\beta\prod_{1\le i<j\le n}{A_{i,j}^{a_{i,j}}}\in \overline{\sigma}^{-1}(\overline{\sigma}(\beta))$, where $\prod_{1\le i<j\le n}{A_{i,j}^{a_{i,j}}}$ belongs to $P_n/[P_n,P_n]$ and $a_{i,j}\in\Z$ for all $1\le i<j\le n$. Let $\theta=\overline{\sigma}(\beta^{-1})$. From Theorem~\ref{teo1} $\alpha_1$ is conjugate to an special braid $\alpha_2=\beta\prod_{A_{r,s}\in T_{\theta}}{A_{r,s}^{S_{r,s}}}$ in $B_n/[P_n,P_n]$, where $T_{\theta}$ is a transversal of the action of $\langle \theta \rangle$ on the set $\mathcal{B}=\{A_{i,j}\mid 1\le i<j\le n\}$, $\mathcal{O}_{\theta}(A_{r,s})$ is the orbit of the element $A_{r,s}$ under this action and $S_{r,s}=\sum_{A_{i,j}\in\mathcal{O}_{\theta}(A_{r,s})}{a_{i,j}}$, for all $A_{r,s}\in T_{\theta}$.
\item[Step 3.]
Finally, using Theorem~\ref{teo2}, we are able to decide if $\alpha_2$ is conjugate to $\beta$, solving a system of equations over the integers.  
As part of this procedure we need to choose $\tilde{C}\in\overline{\sigma}^{-1}(C_{S_n}(\theta))$, where $C_{S_n}(\theta)$ is the centralizer of $\theta$ in $S_n$, and to write $\tilde{C}\beta\tilde{C}^{-1}$ in the form $\beta\prod_{1\le i<j\le n}{A_{i,j}^{c_{i,j}}}$. 
\end{itemize}

Let $\beta\in B_n/[P_n,P_n]$ such that $\overline{\sigma}(\beta^{-1})=\theta=(1,2,\dots,n)$. 
It is not difficult to show that the elements of $C_{S_n}(\theta)$ are $\theta^t$ with $1\leq t\leq n$. For any $t\in \{1,\ldots, n\}$ let $\beta'$ be a coset of a braid such that $\overline{\sigma}(\beta')=\theta^t$, then up to conjugacy, $\beta'=\beta\prod_{1\le i<j\le n}{A_{i,j}^{a_{i,j}}}$.
Since the set $\overline{\sigma}^{-1}(C_{S_n}(\theta))$ is well understood 
in this case, from Theorems~\ref{teo1} and~\ref{teo2}, we know completely the elements of the conjugacy class of $\beta\prod_{A_{t,q}\in T_{\theta}}{A_{t,q}^{z_{t,q}}}$.

\begin{Cor}
Let $n\ge 3$, and let $\beta\in B_n/[P_n,P_n]$ such that $\overline{\sigma}(\beta^{-1})=\theta=(1,2,\dots,n)\in S_n$. The elements $\beta\prod_{A_{t,q}\in T_{\theta}}{A_{t,q}^{z_{t,q}}}$ and $\beta\prod_{A_{t,q}\in T_{\theta}}{A_{t,q}^{b_{t,q}}}$, where $z_{t,q},b_{t,q}\in\Z$ for all $A_{t,q}\in T_{\theta}$, are conjugate if and only if $z_{t,q}=b_{t,q}$ for $A_{t,q}\in T_{\theta}$. 
\end{Cor}
\begin{proof}
We observe that in this case the equality $C_{B_n/[P_n,P_n]}(\beta)=\overline{\sigma}^{-1}(C_{S_n}(\theta))$ holds. For $\tilde{C}\in\overline{\sigma}^{-1}(C_{S_n}(\theta))$, we have $\tilde{C}=\beta^k$ for some $k\in\Z$, thus $A_{C^{-1}(t),C^{-1}(q)}\in\mathcal{O}_{\theta}(A_{t,q})$ for all $A_{t,q}\in T_{\theta}$, where $\overline{\sigma}(\tilde{C})=C$. 
\end{proof}

Next, we illustrate with an example how Theorems~\ref{teo1} and~\ref{teo2} may be applied to determine whether two braids in the Artin braid group quotient $B_n/[P_n,P_n]$ are conjugate.

\begin{Ex}
Are the elements $\sigma_1A_{1,2}$ and $\sigma_1\sigma_2\sigma_1A_{1,3}$ conjugate in $B_3/[P_3,P_3]$? 
To answer this question first we note that $\overline{\sigma}(\sigma_1\sigma_2\sigma_1A_{1,3})=(1,3)$ and $(2,3)(1,3)(2,3)=(1,2)$, so $\overline{\sigma}(\sigma_2\sigma_1\sigma_2\sigma_1A_{1,3}\sigma_2^{-1})=(1,2)$. 
Since $\sigma_2\sigma_1\sigma_2\sigma_1A_{1,3}\sigma_2^{-1}$ and $\sigma_1\sigma_2\sigma_1A_{1,3}$ are conjugate it is sufficient to answer if $\sigma_2\sigma_1\sigma_2\sigma_1A_{1,3}\sigma_2^{-1}$ and $\sigma_1A_{1,2}$ are conjugate in $B_3/[P_3,P_3]$. 
Notice that
\begin{align*}
\sigma_2\sigma_1\sigma_2\sigma_1A_{1,3}\sigma_2^{-1}&=\sigma_1\sigma_2\sigma_1^2A_{1,3}\sigma_2^{-1}\\
&=\sigma_1A_{1,2}A_{1,3}.
\end{align*}
Thus, we must answer if the elements $\sigma_1A_{1,2}A_{1,3}$ and $\sigma_1A_{1,2}$ are conjugate. We will denote by $\theta=\overline{\sigma}(\sigma_1)^{-1}$ and $T_{\theta}=\{A_{1,2},A_{1,3}\}$. 
Notice that $\sigma_1\in\overline{\sigma}^{-1}(C_{S_3}(\theta))$ and $\sigma_1\sigma_1A_{1,2}\sigma_{1}^{-1}=\sigma_1A_{1,2}$. By Theorem~\ref{teo2} $(ii)$ for $\sigma_1A_{1,2}^{b_{1,2}}A_{1,3}^{b_{1,3}}$ to be conjugate to $\sigma_1A_{1,2}$ we should have $1=b_{1,2}$ and $0=b_{1,3}$. Thus, $\sigma_1A_{1,2}$ and $\sigma_1\sigma_2\sigma_1A_{1,3}$ are not conjugate in $B_3/[P_3,P_3]$.
\end{Ex}


\begin{Obs}
Using the ideas from Theorems~\ref{teo1} and~\ref{teo2} to determine whether two elements of $B_n/[P_n,P_n]$ are conjugate, illustrated in the previous two examples, we partition the group $B_3/[P_3,P_3]$ by fibers, coming from the conjugacy classes of elements in $S_3$, such that each fiber was partitioned by conjugacy classes. See Figure~\ref{feijão}.
\begin{figure}[h]
\begin{center}
\includegraphics[scale=0.3]{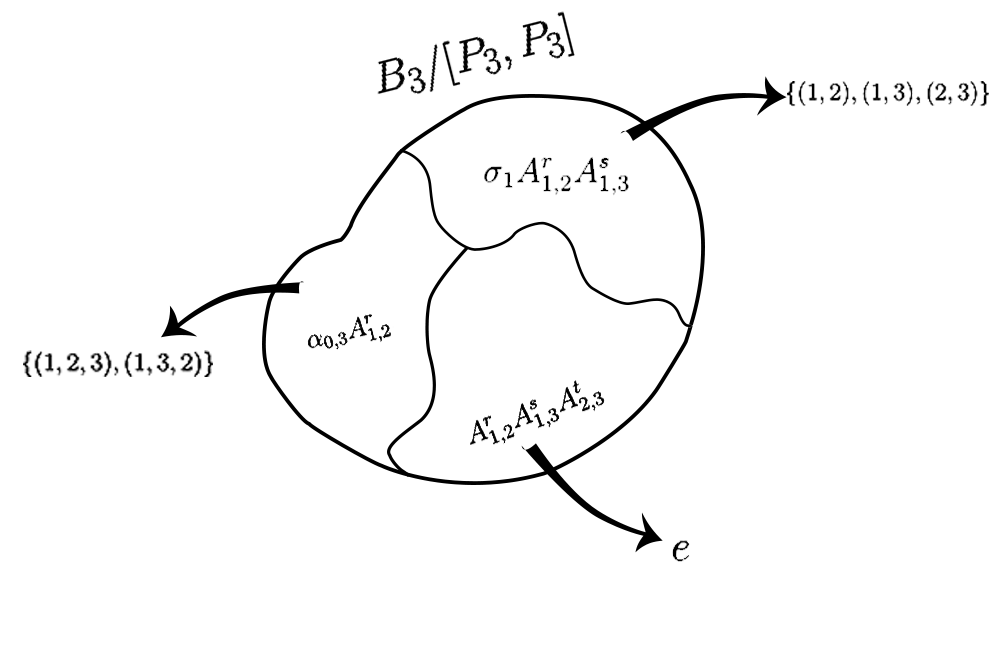}
\caption{Conjugacy classes of elements in the group $B_3/[P_3,P_3]$}
\label{feijão}
\end{center}
\end{figure}\\
In Figure~\ref{feijão} for each $r\in\Z$, the element $\alpha_{0,3}A_{1,2}^r$ represents a conjugacy class of the elements that are projected in the 3-cycle permutations, where $\alpha_{0,3}=\sigma_1\sigma_2$. In its turn, for each $r,s\in\Z$ the element $\sigma_1A_{1,2}^rA_{1,3}^s$ represents a conjugacy class of elements that are projected in the transpositions of $S_3$. Finally, for $r,s,t\in\Z$ such that $r<s<t$ we have the representatives of conjugacy classes of elements that are projected in the trivial element $e$, given of  the form $A_{1,2}^rA_{1,3}^sA_{2,3}^t$.
\end{Obs}




\section{Virtually cyclic subgroups in $B_n/[P_n,P_n]$}\label{sec:3}

Realization of subgroups in the quotient $B_n/[P_n,P_n]$ is not an easy task. Embeddings of finite abelian groups of odd order in $B_n/[P_n,P_n]$ were given in \cite{art} as well as an explicit embedding of a non-abelian group of order 21 (that belongs to a family of groups called \emph{Frobenius group}). After that, many results of \cite{art} were extended to braid groups of complex reflection groups (a generalization of braid groups), in particular the existence of finite subgroups in quotients of these generalized braid groups, see  \cite{bm} and \cite{ma}. 
In \cite{art2} the authors used tools from cohomology theory in order to show that it is possible to realize non-abelian finite groups of odd order in certain quotients of Artin braid groups, in particular in $B_n/[P_n,P_n]$, and some explicit embeddings were showed. 
The following definition can be found in the introduction of \cite{fajo}.


\begin{Def}
A group $G$ is virtually cyclic if it has a cyclic subgroup $H$ of finite index.
\end{Def}

In this section we deal with the problem to realize explicitly infinite virtually cyclic subgroups in the crystallographic group $B_n/[P_n,P_n]$.
Since every finite group is virtually cyclic, in this paper we are interested in infinite virtually cyclic groups. 
There are two distinct types of infinite virtually cyclic groups, the semi-direct products and the amalgamated products as the following result shows. 
\begin{Teo}[{\cite[Theorem~17]{virtual}}]\label{teo17}
Let $G$ be a group. Then the following statements are equivalent.
\begin{itemize}
	\item[(a)] $G$ is a group with two ends.
	\item[(b)] $G$ is an infinite virtually cyclic group.
	\item[(c)] $G$ has a finite normal subgroup $F$ such that $G/F$ is isomorphic to $\Z$ or to the infinite dihedral group $\Z_2\ast \Z_2$.
\end{itemize}
Equivalently, $G$ is of the form:
\begin{itemize}
	\item[(i)] $F\rtimes_{\theta}\Z$ for some action $\theta\in Hom(\Z,\, Aut(F))$, or
	\item[(ii)] $G_1\ast_{F}G_2$ where $[G_i:F]=2$, for $i=1,2$, where $G_1$, $G_2$ and $F$ are finite groups.
\end{itemize}
\end{Teo}
We note that Theorem~\ref{teo17} is well known, most of the first part is due to Epstein \cite{e} and Wall \cite{w}, the proof of the implication (c) implies (b) and of the second part can be found in \cite{virtual}.

\begin{Obs}
Let $n\geq 3$. Since $B_n/[P_n, P_n]$ has no 2-torsion (see \cite[Theorem~2]{art}) the group $B_n/[P_n, P_n]$ has no virtually cyclic subgroups of type amalgamated products. 
Hence, we have the problem to realize infinite virtually cyclic subgroups of $B_n/[P_n, P_n]$ that are semi-direct products.
\end{Obs}

\subsection{Subgroups $\Z_p\rtimes_{\alpha}\Z$ in $B_n/[P_n,P_n]$}

Let $p$ be an odd prime number. Given the homomorphism $\alpha\colon\Z\longrightarrow\Aut(\Z_p)$, note that $\alpha$ is completely determined by the image of the generator of $\Z$. We will denote by $\alpha_k$ the elements of $\Aut(\Z_p)$, where $\alpha_k\colon\Z_p\longrightarrow\Z_p,\lambda\mapsto\lambda^k$ for $k=1,\dots,p-1$ and by $\Z_p\rtimes_{\alpha_k}\Z$ the group $\Z_p\rtimes_{\alpha}\Z$ when $\alpha(x)=\alpha_k$ for $\Z=\langle x\rangle$. 

If $p,n\ge 3$ and $r\ge 0$ are integers such that $p$ is odd and $r+p\le n$, we define $\delta_{r,p},\alpha_{r,p}\in B_n/[P_n,P_n]$ by
\begin{equation}\label{delta}
\delta_{r,p}=\sigma_{r+p-1}\cdots\sigma_{r+\frac{p+1}{2}}\sigma^{-1}_{r+\frac{p+1}{2}}\cdots\sigma^{-1}_{r+1}\mbox{ and } \alpha_{r,p}=\sigma_{r+1}\cdots\sigma_{r+p-1}.
\end{equation}
The elements $\delta_{r,p},\alpha_{r,p}$ were defined in \cite[p.~410]{art} to study the existence of finite order elements in $B_n/[P_n,P_n]$ and a characterization of them.


\begin{proof}[Proof of Theorem~\ref{teoprimos}]
We shall prove that the group $\Z_p\rtimes_{\alpha_k}\Z$ can be realized in $B_n/[P_n,P_n]$ where $k=1,\dots,p-1$. We use the following presentation of the virtually cyclic group $\Z_p\rtimes_{\alpha_k}\Z=\langle A,B\mid A^p=1,BAB^{-1}=A^{k}\rangle$. 
Since $\delta_{0,p}$ (defined in~(\ref{delta}))  have order $p$ in $B_n/[P_n,P_n]$ (see \cite[Lemma~28]{art}), we consider $A=\delta_{0,p}$. Now, we shall prove that there is $\beta\in B_n/[P_n,P_n]$ such that $\Z=\langle\beta\rangle$ and $\beta\delta_{0,p}\beta^{-1}=\delta_{0,p}^{k}$. 
In case there is such element, we should have $\overline{\sigma}(\beta)\overline{\sigma}(\delta_{0,p})\overline{\sigma}(\beta)^{-1}=\overline{\sigma}(\delta_{0,p})^{k}$. 
Notice that $\overline{\sigma}(\delta_{0,p})=(1,\dots,p)$. Since $(1,\dots,p)$ and $(1,\dots,p)^{k}$ are conjugate in $S_n$, there is $\overline{\gamma}\in S_n$ such that $\overline{\gamma}(1,\dots,p)\overline{\gamma}^{-1}=(1,\dots,p)^{k}$. Let $\gamma\in B_n/[P_n,P_n]$ such that $\overline{\sigma}(\gamma)=\overline{\gamma}$. Consider $\beta=\left(\prod_{1\le i<j\le n}{A_{i,j}^{x_{i,j}}}\right)\gamma$ where $x_{i,j}\in\Z$ for all $1\le i<j\le n$, so
\begin{align*}
\beta\delta_{0,p}\beta^{-1}=\delta_{0,p}^{k}&\iff\left(\prod_{1\le i<j\le n}{A_{i,j}^{x_{i,j}}}\right)\gamma\delta_{0,p}\gamma^{-1}\left(\prod_{1\le i<j\le n}{A_{i,j}^{-x_{i,j}}}\right)\delta_{0,p}^{-k}=1.
\end{align*}
Since the homomorphism $\overline{\sigma}\colon B_n/[P_n,P_n]\to S_n$ is surjective with kernel $P_n/[P_n,P_n]$ (see equation~\eqref{sqpura}), $\overline{\sigma}(\gamma\delta_{0,p}\gamma^{-1})=(1,\dots,p)^{k}$ and $\overline{\sigma}(\delta_{0,p})=(1,\dots,p)$, we conclude that $\gamma\delta_{0,p}\gamma^{-1}=\left( \prod_{1\le i<j\le n}{A_{i,j}^{c_{i,j}}}\right)\delta_{0,p}^{k}$, for some $c_{i,j}\in\Z$. 
Thus, applying equation~\eqref{eq:1.2}, we have
\begin{align*}
\beta\delta_{0,p}\beta^{-1}=\delta_{0,p}^{k}&\iff\prod_{1\le i<j\le n}{A_{i,j}^{x_{i,j}}}\left(\prod_{1\le i<j\le n}{A_{i,j}^{c_{i,j}}}\right)\delta_{0,p}^{k}\left(\prod_{1\le i<j\le n}{A_{i,j}^{-x_{i,j}}}\right)\delta_{0,p}^{-k}=1\\
&\iff\prod_{1\le i<j\le n}{A_{i,j}^{x_{i,j}}}\prod_{1\le i<j\le n}{A_{i,j}^{c_{i,j}}}\prod_{1\le i<j\le n}{A_{i,j}^{-x_{\theta(i),\theta(j)}}}=1,
\end{align*} 
where $\theta=\overline{\sigma}(\delta_{0,p}^k)$. The last equation is equivalent to the following system of equations of integers, because $P_n/[P_n,P_n]$ is a torsion free abelian group with basis $\{A_{i,j}\}_{1\leq i<j\leq n}$
\begin{equation}\label{sis1}
x_{i,j}+c_{i,j}-x_{\theta(i),\theta(j)}=0,\mbox{ for all }1\le i<j\le n.
\end{equation}
We can rewrite the system of equations~(\ref{sis1}) as follows. Given $A_{t,r}\in T_{\theta}$ we have the subsystem of equations of integers
$$
\left\{\begin{array}{lc}
x_{t,r}+c_{t,r}-x_{\theta(t),\theta(r)}=0\\
x_{\theta^{k}(t),\theta^{k}(r)}+c_{\theta^{k}(t),\theta^{k}(r)}-x_{\theta^{k+1}(t),\theta^{k+1}(r)}=0,
\end{array}\right.
$$
where $1\le k\le |\mathcal{O}_\theta(A_{t,r})|-1$. We note that each subsystem of equations has solution if and only if $c_{t,r}+c_{\theta(t),\theta(r)}+\cdots+c_{\theta^l(t),\theta^l(r)}=0$ for $l=|\mathcal{O}_\theta(A_{t,r})|-1$. Since $\delta_{0,p}$ have finite order we have $c_{t,r}+c_{\theta(t),\theta(r)}+\cdots+c_{\theta^l(t),\theta^l(r)}=0$ for all $A_{t,r}\in T_{\theta}$ (see \cite[Proposition~29]{art}). Hence the system of equations given in~(\ref{sis1}) have solution. Therefore there is $\beta\in B_n/[P_n,P_n]$ such that $\Z=\langle\beta\rangle$ and $\beta\delta_{0,p}\beta^{-1}=\delta_{0,p}^{k}$.
\end{proof}


\begin{Obs}\label{passos}
Taking $x_{\theta^{k+1}(t),\theta^{k+1}(r)}=c_{\theta^{k}(t),\theta^{k}(r)}+x_{\theta^{k}(t),\theta^{k}(r)}$ for $0\le k\le |\mathcal{O}_\theta(A_{t,r})|-2$ and $x_{t,r}\in\Z$ for all $A_{t,r}\in T_{\theta}$, we have a solution for the system of equations in~$(\ref{sis1})$.

We note that to find explicit generators for the group $\Z_p\rtimes_{\alpha_k}\Z\leq B_n/[P_n,P_n]$ we may exhibit $\gamma\in B_n/[P_n,P_n]$ as in the proof of Theorem~\ref{teoprimos}. In order to do this write $\gamma\delta_{0,p}\gamma^{-1}\delta_{0,p}^{-k}$ as a product of elements of $P_n/[P_n,P_n]$, using \cite[Proposition~15]{art}, and finally solve the system of equations~(\ref{sis1}).
\end{Obs}


\subsection{The virtually cyclic subgroups of $B_n/[P_n,P_n]$ for small values of $n$}

Follows from \cite[Theorems~2~and~6]{art} that every finite subgroup of $B_n/[P_n,P_n]$, for $n=3,4,5$ and $6$, is isomorphic to an abelian group of odd order that embeds in $S_n$. 

In this subsection we exhibit a presentation for all infinite virtually cyclic subgroups in $B_n/[P_n,P_n]$ for $n=3,4,5$ and $6$, up to isomorphism. In order to do this we will follow the steps of Remark~\ref{passos}. We shall illustrate the procedure with an example and note that the other computations are given in a similar way.
We register the information of the presentation of the virtually cyclic subgroups of $B_n/[P_n,P_n]$ for small values of $n$ in tables at the end of this section.

\begin{Obs} We notice that the element $\sigma_5\sigma_4^{-1}\in B_6/[P_6,P_6]$ have finite order. To see this, let $\gamma=\sigma_3\sigma_2\sigma_4\sigma_1\sigma_3\sigma_5\sigma_2\sigma_4\sigma_3$ and note that $\overline{\sigma}(\gamma\sigma_5\sigma_4^{-1}\gamma^{-1})=\overline{\sigma}(\delta_{0,3})$. Thus, $\gamma\sigma_5\sigma_4^{-1}\gamma^{-1}\delta_{0,3}^{-1}\in P_6/[P_6,P_6]$. Computing the crossing numbers of the class representative of the coset of this pure braid 
and applying \cite[Proposition~15]{art} we may conclude that $\gamma\sigma_5\sigma_4^{-1}\gamma^{-1}\delta_{0,3}^{-1}=1$. Since $\delta_{0,3}$ have order $3$ in $B_6/[P_6,P_6]$ (see \cite[Lemma~28]{art}), the element $\sigma_5\sigma_4^{-1}\in B_6/[P_6,P_6]$ also have order $3$.  
\end{Obs}

\begin{Ex}\label{exmvcg}
We exhibit a presentation of the group $(\Z_3\times\Z_3)\rtimes_{\upsilon_2}\Z$ in $B_6/[P_6,P_6]$, where $\upsilon_{2}\colon\Z_3\times\Z_3\longrightarrow\Z_3\times\Z_3$, $\delta_{0,3}\mapsto\sigma_5\sigma_4^{-1},\sigma_5\sigma_4^{-1}\mapsto\delta_{0,3}$. We prove that there is $\beta\in B_6/[P_6,P_6]$ such that $\beta\delta_{0,3}\beta^{-1}=\sigma_5\sigma_4^{-1}$ and $\beta\sigma_5\sigma_4^{-1}\beta^{-1}=\delta_{0,3}$. Let $\gamma=\sigma_3\sigma_2\sigma_4\sigma_1\sigma_3\sigma_5\sigma_4\sigma_2\sigma_3$ and note that $\overline{\sigma}(\gamma\delta_{0,3}\gamma^{-1})=\overline{\sigma}(\sigma_5\sigma_4^{-1})$ and $\overline{\sigma}(\gamma\sigma_5\sigma_4^{-1}\gamma^{-1})=\overline{\sigma}(\delta_{0,3})$. Computing the crossing numbers of the braids $\gamma\delta_{0,3}\gamma^{-1}(\sigma_5\sigma_4^{-1})^{-1},\gamma\sigma_5\sigma_4^{-1}\gamma^{-1}\delta_{0,3}^{-1}\in P_6/[P_6,P_6]$ and applying \cite[Proposition~15]{art} we may conclude that $\gamma\delta_{0,3}\gamma^{-1}(\sigma_5\sigma_4^{-1})^{-1}=1$ and $\gamma\sigma_5\sigma_4^{-1}\gamma^{-1}\delta_{0,3}^{-1}=1$ in $P_6/[P_6,P_6]$. 
Consider the element $\beta=\left(\prod_{1\le i<j\le 6}{A_{i,j}^{x_{i,j}}}\right)\gamma$, for $x_{i,j}\in\Z$. 
Since $\gamma\delta_{0,3}\gamma^{-1}=\sigma_5\sigma_4^{-1}$ in $B_6/[P_6,P_6]$, then
\begin{align*}
\beta\delta_{0,3}\beta^{-1}=\sigma_5\sigma_4^{-1}&\iff\left(\prod_{1\le i<j\le 6}{A_{i,j}^{x_{i,j}}}\right)\gamma\delta_{0,3}\gamma^{-1}\left(\prod_{1\le i<j\le 6}{A_{i,j}^{-x_{i,j}}}\right)(\sigma_5\sigma_4^{-1})^{-1}=1\\
&\iff\left(\prod_{1\le i<j\le 6}{A_{i,j}^{x_{i,j}}}\right)\sigma_5\sigma_4^{-1}\left(\prod_{1\le i<j\le 6}{A_{i,j}^{-x_{i,j}}}\right)(\sigma_5\sigma_4^{-1})^{-1}=1.
\end{align*}
Last equation is equivalent to the system of equations
\begin{eqnarray*}\label{tip4}
\left\{\begin{array}{lc}
x_{1,2}-x_{1,2}=0\\
x_{1,3}-x_{1,3}=0\\
x_{2,3}-x_{2,3}=0
\end{array}\right.,\,
\left\{\begin{array}{lc}
x_{1,4}-x_{1,5}=0\\
x_{1,5}-x_{1,6}=0\\
x_{1,6}-x_{1,4}=0\\
x_{2,4}-x_{2,5}=0\\
x_{2,5}-x_{2,6}=0\\
x_{2,6}-x_{2,4}=0\\
x_{3,4}-x_{3,5}=0\\
x_{3,5}-x_{3,6}=0\\
x_{3,6}-x_{3,4}=0
\end{array}\right.,\,
\left\{\begin{array}{lc}
x_{4,5}-x_{5,6}=0\\
x_{5,6}-x_{4,6}=0\\
x_{4,6}-x_{4,5}=0
\end{array}\right.
\end{eqnarray*}

Now, since $\gamma\sigma_5\sigma_4^{-1}\gamma^{-1}=\delta_{0,3}$ in $B_6/[P_6, P_6]$ as mentioned above in this example, we obtain the following  
\begin{align*}
\beta\sigma_5\sigma_4^{-1}\beta^{-1}=\delta_{0,3}&\iff\left(\prod_{1\le i<j\le 6}{A_{i,j}^{x_{i,j}}}\right)\gamma\sigma_5\sigma_4^{-1}\gamma^{-1}\left(\prod_{1\le i<j\le 6}{A_{i,j}^{-x_{i,j}}}\right)\delta_{0,3}^{-1}=1\\
&\iff\left(\prod_{1\le i<j\le 6}{A_{i,j}^{x_{i,j}}}\right)\delta_{0,3}\left(\prod_{1\le i<j\le 6}{A_{i,j}^{-x_{i,j}}}\right)\delta_{0,3}^{-1}=1.
\end{align*}
Last equation is equivalent to the system of equations
\begin{eqnarray*}\label{tip5}
\left\{\begin{array}{lc}
x_{1,2}-x_{2,3}=0\\
x_{2,3}-x_{1,3}=0\\
x_{1,3}-x_{1,2}=0\\
\end{array}\right.,\,
\left\{\begin{array}{lc}
x_{1,4}-x_{2,4}=0\\
x_{2,4}-x_{3,4}=0\\
x_{3,4}-x_{1,4}=0\\
x_{1,5}-x_{2,5}=0\\
x_{2,5}-x_{3,5}=0\\
x_{3,5}-x_{1,5}=0\\
x_{1,6}-x_{2,6}=0\\
x_{2,6}-x_{3,6}=0\\
x_{3,6}-x_{1,6}=0
\end{array}\right.,\,
\left\{\begin{array}{lc}
x_{4,5}-x_{4,5}=0\\
x_{4,6}-x_{4,6}=0\\
x_{4,6}-x_{4,6}=0
\end{array}\right.
\end{eqnarray*}
Therefore, taking $x_{1,2}=1,x_{1,4}=1$ and $x_{4,5}=1$ we obtain 
$$
(\Z_3\times\Z_3)\rtimes_{\alpha_2}\Z\cong\langle\delta_{0,3},\sigma_5\sigma_4^{-1},(A_{1,2}A_{1,3}A_{2,3})(A_{4,5}A_{4,6}A_{5,6})(A_{2,4}A_{2,5}A_{2,6}A_{3,4}A_{3,5}A_{3,6}A_{1,4}A_{1,5}A_{1,6})\gamma\rangle.
$$
\end{Ex}

\begin{Obs}
Note that the order of $\Aut(\Z_3\times\Z_3)$ is six. But it is not possible  to realize the groups $(\Z_3\times\Z_3)\rtimes_{\upsilon_5}\Z$ and $(\Z_3\times\Z_3)\rtimes_{\upsilon_6}\Z$ in $B_6/[P_6,P_6]$, where $\upsilon_5\colon\Z_3\times\Z_3\longrightarrow\Z_3\times\Z_3,\delta_{0,3}\mapsto\delta_{0,3}(\sigma_5\sigma_4^{-1})^2,\sigma_5\sigma_4^{-1}\mapsto\delta_{0,3}^2\sigma_5\sigma_4^{-1}$ and $\upsilon_6\colon\Z_3\times\Z_3\longrightarrow\Z_3\times\Z_3,\delta_{0,3}\mapsto\delta_{0,3}^2\sigma_5\sigma_4^{-1},\sigma_5\sigma_4^{-1}\mapsto\delta_{0,3}(\sigma_5\sigma_4^{-1})^2$, since in each case the system obtained has no solution in $\Z$. 
\end{Obs}
The other cases are presented in the following tables. The computations are similar to the ones showed in Example~\ref{exmvcg} in order to obtain the generators of the virtually cyclic group. The element $\gamma\in B_n/[P_n,P_n]$ given in the tables is the one that appears in the proof of Theorem~\ref{teoprimos} (see also Remark~\ref{passos}).

\begin{table}[htb]
\caption{The virtually cyclic subgroups in $B_3/[P_3,P_3]$}
\label{tab3}
\large
\centering
{\tiny
\begin{tabular}{|c|c|c|c|}
\hline
Group & Action & Generators & $\gamma$ \\
\hline
\hline
$\Z_3\rtimes\Z $&$\alpha_{1}\colon\Z_3\longrightarrow \Z_3$, $\delta_{0,3}\mapsto\delta_{0,3}$&$\delta_{0,3},A_{1,2}A_{2,3}A_{1,3}$ & $1$ \\
\cline{2-4}
&$\alpha_{2}\colon\Z_3\longrightarrow \Z_3$, $\delta_{0,3}\mapsto\delta_{0,3}^2$& $\delta_{0,3},A_{1,2}A_{1,3}A_{2,3}\gamma$ & $\sigma_1\sigma_2\sigma_1$ \\
\hline
\end{tabular}
}
\end{table}
\begin{table}[htb]
\caption{The virtually cyclic subgroups in $B_4/[P_4,P_4]$}
\label{tab3}
\large
\centering
{\tiny
\begin{tabular}{|c|c|c|c|}
\hline
Group & Action & Generators & $\gamma$ \\
\hline
\hline
$\Z_3\rtimes\Z $&$\alpha_{1}\colon\Z_3\longrightarrow \Z_3$, $\delta_{0,3}\mapsto\delta_{0,3}$ &$\delta_{0,3},(A_{1,2}A_{2,3}A_{1,3})(A_{1,4}A_{2,4}A_{3,4})$& $1$ \\
\cline{2-4}
&$\alpha_{2}\colon\Z_3\longrightarrow \Z_3$, $\delta_{0,3}\mapsto\delta_{0,3}^2$&$\delta_{0,3},(A_{1,2}A_{1,3}A_{2,3})(A_{1,4}A_{3,4}A_{2,4})\gamma$ &$\sigma_1\sigma_2\sigma_1$  \\
\hline
\end{tabular}
}
\end{table}
\begin{table}[htb]
\caption{The virtually cyclic subgroups in $B_5/[P_5,P_5]$ }
\label{tab2}
\large
\centering
{\tiny
\begin{tabular}{|c|c|c|c|}
\hline
Group & Action & Generators & $\gamma$ \\
\hline
\hline
$\Z_3\rtimes\Z $&$\alpha_{1}\colon\Z_3\longrightarrow\Z_3$, $\delta_{0,3}\mapsto\delta_{0,3}$ &$\delta_{0,3},(A_{1,2}A_{2,3}A_{1,3})(A_{1,4}A_{2,4}A_{3,4})(A_{1,5}A_{2,5}A_{3,5})(A_{4,5})$ & $1$\\
\cline{2-4}
&$\alpha_{2}\colon\Z_3\longrightarrow \Z_3$, $\delta_{0,3}\mapsto\delta_{0,3}^2$& $\delta_{0,3},(A_{1,2}A_{1,3}A_{2,3})(A_{1,4}A_{3,4}A_{2,4})(A_{1,5}A_{3,5}A_{2,5})(A_{4,5})\gamma$ & $\sigma_1\sigma_2\sigma_1$ \\
\hline
$\Z_5\rtimes\Z$& $\alpha_{1}\colon\Z_5\longrightarrow \Z_5$, $\delta_{0,5}\mapsto\delta_{0,5}$&$\delta_{0,5},(A_{1,2}A_{2,3}A_{3,4}A_{4,5}A_{1,5})(A_{2,4}A_{3,5}A_{1,4}A_{2,5}A_{1,3})$ & $1$ \\
\cline{2-4}
&$\alpha_{2}\colon\Z_5\longrightarrow \Z_5$, $\delta_{0,5}\mapsto\delta_{0,5}^2$&$\delta_{0,5},(A_{1,2}A_{3,4})(A_{2,4}A_{1,4}A_{1,3}A_{3,5})\gamma$ &$\sigma_4\sigma_2\sigma_3$   \\
\cline{2-4}
&$\alpha_{3}\colon\Z_5\longrightarrow \Z_5$, $\delta_{0,5}\mapsto\delta_{0,5}^3$&$\delta_{0,5},(A_{1,2}A_{4,5}A_{2,3}A_{1,5}A_{3,4})(A_{2,4}A_{2,5}A_{3,5}A_{1,3}^2A_{1,4}A_{2,4})\gamma$ & $\sigma_3\sigma_4\sigma_2$ \\
\cline{2-4}
&$\alpha_{4}\colon\Z_5\longrightarrow \Z_5$, $\delta_{0,5}\mapsto\delta_{0,5}^4$&$\delta_{0,5},(A_{1,2})(A_{2,4}A_{1,3}^2A_{2,5}A_{1,4}A_{3,5})\gamma$ & $\sigma_4\sigma_3\sigma_2\sigma_3\sigma_4\sigma_3$ \\
\hline
\end{tabular}
}
\end{table}
\begin{table}[htb]
\caption{The virtually cyclic subgroups in $B_6/[P_6,P_6]$ }
\label{tab1}
\large
\centering
{\tiny
\begin{tabular}{|c|c|}
\hline
Group & Action \\
\hline
\hline
$\Z_3\rtimes\Z$ &$\alpha_{1}\colon\Z_3\longrightarrow \Z_3,\delta_{0,3}\mapsto\delta_{0,3}$\\
\cline{2-2}
&$\alpha_{2}\colon\Z_3\longrightarrow \Z_3$, $\delta_{0,3}\mapsto\delta_{0,3}^2$  \\
\hline
$\Z_3\rtimes\Z $&$\alpha_{1}\colon\Z_3\longrightarrow\Z_3$, $\delta_{0,3}\sigma_5\sigma_4^{-1}\mapsto\delta_{0,3}\sigma_5\sigma_4^{-1}$ \\
\cline{2-2}
&$\alpha_{2}\colon\Z_3\longrightarrow \Z_3$, $\delta_{0,3}\sigma_5\sigma_4^{-1}\mapsto\delta_{0,3}^2(\sigma_5\sigma_4^{-1})^2$ \\
\hline
$\Z_5\rtimes\Z$ &$\alpha_{1}\colon\Z_5\longrightarrow \Z_5$,$\delta_{0,5}\mapsto\delta_{0,5}$ \\
\cline{2-2}
&$\alpha_{2}\colon\Z_5\longrightarrow \Z_5$, $\delta_{0,5}\mapsto\delta_{0,5}^2$ \\
\cline{2-2}
&$\alpha_{3}\colon\Z_5\longrightarrow \Z_5$, $\delta_{0,5}\mapsto\delta_{0,5}^3$ \\
\cline{2-2}
&$\alpha_{4}\colon\Z_5\longrightarrow \Z_5$, $\delta_{0,5}\mapsto\delta_{0,5}^4$  \\
\hline
$(\Z_3\times\Z_3)\rtimes\Z$& $\upsilon_{1}\colon\Z_3\times\Z_3\longrightarrow\Z_3\times\Z_3$, $\delta_{0,3}\mapsto\delta_{0,3},\sigma_5\sigma_4^{-1}\mapsto\sigma_5\sigma_4^{-1}$ \\
\cline{2-2}
&$\upsilon_{2}\colon\Z_3\times\Z_3\longrightarrow\Z_3\times\Z_3$, $\delta_{0,3}\mapsto\sigma_5\sigma_4^{-1},\sigma_5\sigma_4^{-1}\mapsto\delta_{0,3}$ \\
\cline{2-2}
&$\upsilon_{3}\colon\Z_3\times\Z_3\longrightarrow\Z_3\times\Z_3$, $\delta_{0,3}\mapsto\delta_{0,3}^2,\sigma_5\sigma_4^{-1}\mapsto(\sigma_5\sigma_4^{-1})^2$ \\
\cline{2-2}
&$\upsilon_{4}\colon\Z_3\times\Z_3\longrightarrow\Z_3\times\Z_3$, $\delta_{0,3}\mapsto(\sigma_5\sigma_4^{-1})^2,\sigma_5\sigma_4^{-1}\mapsto\delta_{0,3}^2$  \\
\hline
\end{tabular}
\begin{tabular}{|c|c|}
\hline
Generators & $\gamma$ \\
\hline
\hline
$\delta_{0,3},(A_{1,2}A_{2,3}A_{1,3})(A_{1,4}A_{2,4}A_{3,4})(A_{1,5}A_{2,5}A_{3,5})(A_{1,6}A_{2,6}A_{3,6})(A_{4,5})(A_{4,6})(A_{5,6})$ & $1$ \\
\hline
$\delta_{0,3},(A_{1,2}A_{1,3}A_{2,3})(A_{1,4}A_{3,4}A_{2,4})(A_{1,5}A_{3,5}A_{2,5})(A_{1,6}A_{3,6}A_{2,6})(A_{4,5})(A_{4,6})(A_{5,6})\gamma$ &$\sigma_1\sigma_2\sigma_1$  \\
\hline
$\delta_{0,3}\sigma_5\sigma_4^{-1},(A_{1,2}A_{2,3}A_{1,3})(A_{4,5}A_{5,6}A_{4,6})(A_{1,4}A_{2,5}A_{3,6})(A_{1,5}A_{2,6}A_{3,4})(A_{1,6}A_{2,4}A_{3,5})$ & $1$ \\
\hline
$\delta_{0,3}\sigma_5\sigma_4^{-1},(A_{1,2}A_{2,3}A_{1,3})(A_{4,5}A_{5,6}A_{4,6})(A_{1,4}A_{2,5}A_{3,6})(A_{1,5}A_{2,6}A_{3,4})(A_{1,6}A_{2,4}A_{3,5})\gamma$  & $\sigma_1\sigma_2\sigma_1\sigma_4\sigma_5\sigma_4$ \\
\hline
$\delta_{0,5},(A_{1,2}A_{2,3}A_{3,4}A_{4,5}A_{1,5})(A_{2,4}A_{3,5}A_{1,4}A_{2,5}A_{1,3})(A_{1,6}A_{2,6}A_{3,6}A_{4,6}A_{5,6})$  & $1$ \\
\hline
$\delta_{0,5},(A_{1,2}A_{3,4})(A_{2,4}A_{1,4}A_{1,3}A_{3,5})(A_{1,6}A_{3,6}A_{5,6}A_{2,6}A_{4,6})\gamma$ &$\sigma_4\sigma_2\sigma_3$ \\
\hline
$\delta_{0,5},(A_{1,2}A_{3,4})(A_{2,4}A_{1,4}A_{1,3}A_{3,5})(A_{1,6}A_{4,6}A_{2,6}A_{5,6}A_{3,6})\gamma$ &$\sigma_3\sigma_4\sigma_2$ \\
\hline
$\delta_{0,5},(A_{1,2}A_{3,4})(A_{2,4}A_{1,4}A_{1,3}A_{3,5})(A_{1,6}A_{5,6}A_{4,6}A_{3,6}A_{2,6})\gamma$ &$\sigma_4\sigma_3\sigma_2\sigma_3\sigma_4\sigma_3$  \\
\hline
$\delta_{0,3},\sigma_5\sigma_4^{-1},(A_{1,3}A_{2,3}A_{1,3})(A_{4,5}A_{5,6}A_{4,6})(A_{2,4}A_{2,5}A_{2,6}A_{3,4}A_{3,5}A_{3,6}A_{1,4}A_{1,5}A_{1,6})$ & $1$ \\
\hline
$\delta_{0,3},\sigma_5\sigma_4^{-1},(A_{1,3}A_{1,3}A_{2,3})(A_{4,5}A_{4,6}A_{5,6})(A_{2,4}A_{2,5}A_{2,6}A_{3,4}A_{3,5}A_{3,6}A_{1,4}A_{1,5}A_{1,6})\gamma$ & $\sigma_3\sigma_2\sigma_4\sigma_1\sigma_3\sigma_5\sigma_4\sigma_2\sigma_3$ \\
\hline
$\delta_{0,3},\sigma_5\sigma_4^{-1},(A_{1,3}A_{2,3}A_{1,3})(A_{4,5}A_{5,6}A_{4,6})(A_{2,4}A_{2,5}A_{2,6}A_{3,4}A_{3,5}A_{3,6}A_{1,4}A_{1,5}A_{1,6})\gamma$ &$\sigma_1\sigma_2\sigma_1\sigma_4\sigma_5\sigma_4$ \\
\hline
$\delta_{0,3},\sigma_5\sigma_4^{-1},(A_{1,3}A_{2,3}A_{1,3})(A_{4,5}A_{5,6}A_{4,6})(A_{2,4}A_{2,5}A_{2,6}A_{3,4}A_{3,5}A_{3,6}A_{1,4}A_{1,5}A_{1,6})\gamma$ &$\sigma_1\sigma_5\sigma_4\sigma_3\sigma_2\sigma_3\sigma_4\sigma_5\sigma_1\sigma_2\sigma_4\sigma_3\sigma_2\sigma_4\sigma_3$  \\
\hline
\end{tabular}
}
\end{table}



%



\end{document}